
%
%
%
%
%
\magnification=1200
\def\cite{\bf}
\def\theorem{\vskip\baselineskip\bf}
\def\demo{\vskip\baselineskip\par{\bf Proof:}}
\def\qed{{\bf Q.E.D.}}
\def\Cal{\sl}
\def\Atilif{\tilde  A}
\def\Ainf{A}

\def\invB{{{\rm inv}}{\tilde B}}
\def\invA{ { {\rm inv}}\tilde A}

\def\R{\bf R}
\def\Z{\bf Z}
\def\N{\bf N}
\def\T{\bf T}
\def\C{\bf C}

\tolerance=10000
\baselineskip=16.5pt plus .2pt  
\parskip=4.5pt plus1.5pt minus 1pt  
%
%
%
\centerline {\bf A Short Proof that $M_{n}(A)$ is local if $A$ is
local and Fr\'echet}
\vskip\baselineskip
\centerline {\bf Larry B. Schweitzer}
%
%
\vskip\baselineskip
\centerline {{\it Current address:}\ \  Department of
Mathematics, University of Victoria}
\centerline{  Victoria, B.C. Canada V8W 3P4}
\centerline{{\it Email Address:}\ \ lschweit@sol.uvic.ca}
\vskip\baselineskip
\centerline{\bf Abstract}
\par
We give a short and very general proof of the fact that the property
of a dense Fr\'echet subalgebra of a Banach algebra
being local, or closed under the holomorphic
functional calculus
in the Banach algebra, is preserved by tensoring with the $n\times n$
matrix algebra of the complex numbers.
\vskip\baselineskip
\vskip\baselineskip
\centerline{ Introduction }
\vskip\baselineskip
\par
We say that a dense subalgebra $A$ of a Banach algebra $B$,
having the same unit as $B$, is local {\cite{[2]}}, or
closed under the holomorphic functional calculus
in $B$ if for
each $a\in A$ and each function $f$ holomorphic in a neighborhood
of the spectrum of $a$ in $B$, the element $f(a)$ of $B$ lies in $A$.
If $A$ is local in $B$, then
the inclusion map $\iota \colon A \hookrightarrow B$ induces an
isomorphism of $K$-theories $\iota_{*}\colon K_{*}(A)\cong
K_{*}(B)$ {\cite{[5, VI.3]}}, {\cite{[3, Theorem A.2.1]}}.
This is one indication of the importance of local dense subalgebras.
\par
In \S 1, we introduce the notion of spectral invariance, and show
that, if $A$ is a Fr\'echet algebra under some possibly stronger
topology than that of $B$, then $A$ is spectral invariant in $B$ iff
$A$ is local in $B$.  We then show that $A$ is spectral invaraint
in $B$ iff
every simple $A$-module (or algebraically
irreducible representation of $A$) is contained in a $B$-module, and
use this to give a simple proof that the
$n\times n$ matrices $M_{n}(A)$ are spectral invariant in
$M_{n}(B)$ iff $A$ is spectral invariant in $B$.
This result is useful in $K$-theory {\cite{[10, Corollary 7.9]}}.
We conclude with several interesting examples.  This paper
is a portion of my Ph.D. thesis, written at the University
of California at Berkeley, under the supervision of Marc A.
Rieffel.
\vskip\baselineskip
\centerline{\S 1 An Equivalent Condition for Spectral Invariance. }
\par
We show in this section that a dense subalgebra $A$ of a topological
$Q$-algebra $B$ is spectral invariant in $B$ if and only if
all simple $A$-modules are contained in $B$-modules.
\par
{\bf Definition 1.1.\ \ }
By a topological algebra, we mean a topological vector space
 over a topological
ring $k$ with unit, with a $k$-algebra structure for which
 multiplication is separately continuous.
Let $B$ be a topological algebra
over a
topological ring $k$ with unit.  Let
$A$ be a subalgebra.
If $B$ has no unit let $\tilde B $ be $B$ with unit adjoined,
 and let $\Atilif$ be $\Ainf$ with the same unit adjoined.
If $B$ has a unit, and $A$ has the same unit, we let
$\Atilif =A$, $\tilde B =B$.  If $B$ has a unit $1_{B}$ and
$A$ has a different unit, we let $\tilde B = B$ and
$\tilde A$ be the subalgebra of $B$
generated by $A$ and $1_{B}$.
Note that in any case $\tilde B$ is then
a topological algebra, and $\Atilif $ a  subalgebra.
If $A$ is dense in $B$, then $\tilde A$ is dense in $\tilde B$
and $B$ has the same unit as $A$ if $A$ is unital.
We often make the additional
assumption that the group $\invB$  of invertible elements
in $\tilde B$ is open.
If $\invB$ is open, we call $B$ a {\it $Q$-algebra}
(see {\cite{[9]}}).
If, in addition, inversion in $\tilde B$ is continuous,
we say that $B$
is a {\it good} topological algebra (see {\cite{[3]}}).
\par
We say that ${\Ainf}$ is a {\it spectral invariant (SI)}
subalgebra of $B$ if the invertible elements of $\tilde A$
are precisely those elements of $\tilde A$ which are invertible
in $\tilde B$.
Note that $A$ is SI in $B$ if and only if for every $a\in A$,
the spectrum of $a$ is the same in $\tilde A$ and in $\tilde B$.
\par
By a Fr\'echet algebra, we shall mean a Fr\'echet space with an
algebra structure for which multiplication is jointly
continuous {\cite{[13]}}.
We say that a Fr\'echet algebra is $m$-convex if its topology
is given by a family of submultiplicative seminorms.
An $m$-convex Fr\'echet algebra has a holomorphic functional
calculus
{\cite{[10, Lemma 1.3]}}.
If $B$ is an $m$-convex Fr\'echet
algebra and $A$ is a subalgebra, we say that {\it $A$ is
closed
under the holomorphic functional calulus in $B$ }
(or {\it local in } $B$ {\cite{[2]}}) if and only if
for every $a \in \tilde A$ and $f$ holomorphic
in a neighborhood of $spec_{\tilde B}(a)$, the element $f(a)$
of $\tilde B$ lies in $\tilde A$.
\par
We say that $A$ is a {\it Fr\'echet subalgebra} of a topological
algebra $B$ if $A$ is a Fr\'echet algebra, and there exists
an injective and continuous algebra homomorphism (called
the inclusion map)
$\iota \colon A \hookrightarrow
B$.
{\theorem{Lemma 1.2.\ \ }}{\it Let $B$ be an  $m$-convex
Fr\'echet
$Q$-algebra and $\Ainf $ a Fr\'echet subalgebra.
Then $\Ainf$ is closed under the holomorphic
functional calculus in $B$ iff
$\Ainf$ is spectral invariant in $B$.}
{\demo}
For $\Rightarrow$, note that
if $a \in \invB \cap \Atilif$,
then $f(\lambda ) = 1/\lambda $
is holomorphic in a neighborhood of
the closed set $spec_{\tilde B}(a)$ (compare
 {\cite{[2, Proposition 3.1.3]}}).
(We  use $\invB$ open to see that $spec_{\tilde B}(a)$
is closed.)
\par
For $\Leftarrow$
first note that ${\rm inv}  \Atilif = \invB \cap \Atilif$
is open since $\invB$ is open and the inclusion
map from $\Ainf$ to $B$ is continuous.
It follows that $A$ has holomorphic functional calculus
{\cite{[8]}}
{\cite{[13]}} which agrees with the holomorphic functional
calculus
on $B$.  This fact is well known, but I have not been
able to find
it stated in the literature, so I will prove it here
for convenience.
  Since $\invA$ is open and $\Atilif$ is a Fr\'echet
algebra, the inverse map $x \mapsto
x^{-1}$ in ${\rm inv}\Atilif$ is continuous
{\cite{[13, Chapter 7, Proposition 2]}}.
 If $f$ is holomorphic in a neighborhood $N$
of ${spec}_{\tilde B}(a)$, then $\lambda \mapsto
f(\lambda)(\lambda
- a)^{-1}$ is a well defined and continuous map from
any closed curve $\Gamma$ contained
in $N-spec_{\tilde B}(a)$.  From
this one checks that the integral expression
$$f(a) = {1\over 2\pi \imath }\int_{\Gamma} f(\lambda)
(\lambda
- a)^{-1} d\lambda$$
has its partial Riemann sums converging in $\Atilif$.
Thus $f(a)
\in \Atilif$.\qed
\par
We shall often call a
simple module an irreducible representation.
Irreducible will mean algebraically irreducible.
We shall use the terminology module-representation and
submodule-subrepresentation interchangably.
\par
To motivate Theorem 1.4, we consider the
situation where the  algebra $B$ and its dense subalgebra
$A$ are commutative and unital Banach algebras.
Then all simple $A$-modules
are in fact homomorphisms of $\Ainf$ to
$\bf C$.  If $a$ is not invertible in $\Ainf$, one can
easily construct a homomorphism from $\Ainf$ to $\bf C$
that has $a$ in its kernel.  If this homomorphism
extends to a homomorphism from $B$ to $\bf C$, clearly
$a$ cannot be invertible in $B$ (so if all homomorphisms
of $\Ainf$ extend to ones of $B$,
$\Ainf$ must be spectral invariant in $B$).
\par
Conversely, let $h$ be a homomorphism of $A$ to $\C$ with
no extension to a homomorphism of $B$.
Assume that the kernel $M$ of $h$ is not dense in $B$.
Then the closure $\overline M$ is a proper ideal in $B$.
By Zorn's lemma, $\overline M$ is contained in a maximal closed
ideal $N$ in $B$. The
quotient map
$$ B\rightarrow B/{N} \cong \C \qquad
{\rm {(Gelfand-Mazur Theorem)}}$$
gives an extension of the homomorphism $h$ (it extends $h$
since $N\cap A=M$ by the maximality of $M$ as an ideal in $A$).
So since the homomorphism $h$ has no extension,
our assumption that $M$ is not dense in $B$ has
led us to a contradiction.
So $M$ must be dense in $B$.  Since $B$ is a $Q$-algebra,
$M$ must contain an invertible element
of $B$ which is contained in $A$.
Since $M$ is a proper ideal of $A$,
this element is not invertible in $\Ainf$, and thus $\Ainf$ is
not spectral invariant in $B$.
So $A$ is SI in $B$ iff every simple $A$-module is
contained in a $B$-module.
We will generalize these observations in Theorem 1.4.
\par In some literature ({\cite{[3]}}, {\cite{[12]}}),
the following more general setting is considered.
\par
{\bf Definition 1.3.\ \ } Let $B$ be a topological algebra
and let
$\varphi\colon A\rightarrow
B$ be an algebra homomorphism.
We define $\tilde A$ and $\tilde B$ as in Definition 1.1.
We say that $A$ is
{\it spectral invariant} in $B$ if for all
$a\in \tilde A$, $a$
is invertible in $\tilde A$ whenever  $\varphi (a)$ is
invertible in $\tilde B$.
We say that an $A$-module $V$ is {\it contained in}
a $B$-module $W$ if
$V\subseteq W$, $V$ factors to a $\varphi (A)$-module,
and the action of
$B$  on
$W$ is consistent with the action of $\varphi (A)$ on $V$.
\vskip\baselineskip
\par
We prove the following theorem in this setting.  It is similar
to {\cite{[7, Theorem 2.4.16]}}.
{\theorem{Theorem 1.4.\ \ }}{\it Let $B$ be a $Q$-algebra,
and let
$\varphi \colon A \rightarrow B$ be an algebra homomorphism
with dense image.
The following are equivalent:
\par
$(\romannumeral 1 )$ $\Ainf$ is  spectral
invariant in $B$.
\par
$(\romannumeral2)$ For every maximal left ideal $N$ in
$\Atilif$,   we have $\varphi^{-1}
({\overline {\varphi(N)}}^{\tilde B})
 = N$.
\par
$(\romannumeral 3)$ Every simple $ \Ainf $-module $V$
is contained in a $B$-module $W$.
If so desired, $W$ can  be taken
to be a simple $B$-module.
\par
The proof of $(\romannumeral3)\Rightarrow (\romannumeral1)$
does not require $\invB$ open or $\varphi(A)$ dense.
\par
Also, the $B$-module $W$ constructed in our proof
(of $(\romannumeral2 )
\Rightarrow
(\romannumeral3)$) inherits a natural
topology, making it a continuous $B$-module for which the
action of $B$ is continuous (more precisely, the
map $(b, w)\mapsto
bw$ from $B \times W$ to $W$ is separately continuous).
In this topology, $V$ is dense in $W$.
}
{\demo}\par
 $(\romannumeral 1)\Rightarrow (\romannumeral2 )$. Assume $\Ainf$
is SI in $B$.
Let $N$ be a maximal left   ideal in
${\Atilif}$.   Extend $\varphi$ in the natural way to a map
$\varphi \colon  {\tilde A} \rightarrow {\tilde B}$, again
having dense image.
Then $N \cap \invA = \emptyset$ so $\varphi(N) \cap \invB =
\emptyset$ by spectral invariance.
Since $\invB$ is open in $\tilde B$, the closure
${\overline {\varphi(N)}}$ of $\varphi(N)$ in
$\tilde B$ is a proper
subset of $\tilde B$, and is also a left   ideal in $\tilde B$.
Since $\varphi(\Atilif)$ is dense in $\tilde B$,
$\varphi^{-1}({\overline {\varphi(N)}})
 \not=
\Atilif$.
Thus $\varphi^{-1}({\overline {\varphi(N)}})$ is a proper
left  ideal in $\Atilif$.  Since this ideal contains $N$ and $N$ is
maximal, we have $\varphi^{-1}({\overline {\varphi(N)}}) = N$.
(Note that any maximal left ideal of $\tilde A$ therefore contains
the kernel of $\varphi$.)
\vskip\baselineskip
\par
$(\romannumeral2 ) \Rightarrow (\romannumeral3)$.  Let
$V$ be a simple
$\Ainf$-module.
Extend $V$ to a $\Atilif$-module in the natural way.
Let $v\in V$, $v\not= 0$.
Let $N$ be the kernel
of the map $\Atilif \to V$, $a\mapsto av$.  This map is
onto by simplicity, so $V\cong \Atilif /N$ as $\Ainf $
modules.  Also by simplicity, $N$ is a maximal left ideal,
so by $(\romannumeral2 )$ we
have $\varphi^{-1}({\overline {\varphi(N)}}) = N$.
Set $M={\overline {\varphi(N)}}$. Then
the map $\varphi \colon {\tilde A}
\rightarrow {\tilde B}$ gives the inclusion
$$V\cong \Atilif/N \subseteq {\tilde B} /M \eqno(*)$$
where $\tilde B/M$ is a $B$-module.  Take $W={\tilde B}/M$.
This proves
$(\romannumeral2 ) \Rightarrow (\romannumeral3)$.
\par
Clearly $M$ is a closed ideal in $\tilde B$ so $W$
inherits a natural
topology in which $V$ is dense.
If we replace $M$ with a maximal
ideal in $\tilde B$ containing
${\overline {\varphi(N)}}$, we still have
 $\varphi^{-1}(M) = N$ and so $(*)$ continues to
hold. Also $M$ must be closed, so $W={\tilde B}/M$ is a
simple (and continuous) $B$-module extending the action
of $\Ainf$ on $V$ and containing $V$ as a
dense subset.
\vskip\baselineskip
\par
$(\romannumeral3) \Rightarrow (\romannumeral 1)$.
%
Let $a \in \Atilif$,
$\varphi(a)^{-1}\in {\tilde B} - \varphi(\Atilif)$.
We show that $(\romannumeral3)$ is not satisfied.
The element $a$ cannot be left invertible in $\Atilif$
(otherwise we would have $\varphi(a)^{-1}
\in \varphi(\tilde A)$). Hence
$<a>={\Atilif}a$
is a proper left ideal in $\Atilif$.  By Zorn's lemma
there is some
maximal left ideal $N$ in $\Atilif$ containing $<a>$.
View $V=\Atilif/N$ as an $A$-module.
Since
$N$ is maximal,
$V$ is simple.
Note that $a[1]=[a]=0$, while $[1]\not= 0$.
Hence, even if the kernel of $\varphi$ acts trivially
on $V$ (so $V$ factors to an $\varphi({\tilde A})$-module),
 $V$ cannot be contained in any $B$-module
since $\varphi(a)$ is invertible in $\tilde B$.
This proves the contrapositive of
$(\romannumeral 3) \Rightarrow  (\romannumeral1)$.  Note
that we did not use the assumption that $\invB$ is open or that
$\varphi(A)$ is dense for this.\qed
{\theorem{Corollary 1.5.\ \ }}{\it Let $B$
be a C*-algebra, and let
$A$ be a dense Fr\'echet subalgebra.
Then the following are equivalent.
\par
{\bf (1).\ \ } $A$ is closed under the
holomorphic functional
calculus in $B$.
\par
{\bf (2).\ \ } Every  simple (possibly discontinuous)
$A$-module is contained as a dense subspace of an irreducible
*-representation of $B$ on a Hilbert space.}
{\demo}
Condition (1) is equivalent to the spectral
invariance of $A$ in $B$
by Lemma 1.2.
Condition (2) is equivalent to condition $(\romannumeral3)$
of Theorem 1.4 if $B$ is a C*-algebra by
{\cite{[6, Corollary 2.9.6(i)]}}.
\qed
\par
As a point of interest, we remark that if $A$ is
a $Q$-algebra, then every
simple $A$-module $V$ is continuous
for a natural topology.  Namely choose $v\in V$ and let $N$ be
the set of $a \in \tilde A$ such that $av=0$. Then $N$ is closed
by maximality and $V\cong {\tilde A}/N$
gives the topology on $V$.
\par
{\bf Example 1.6.\ \ }  For an example
where $\varphi$ is not injective,
and $A$ is SI in $B$,
let $A$
be power series in two variables, $\C[[X,Y]]$, and
let $B$ be power series in one variable, $\C[[X]]$.  Define
$$ \varphi (\sum c_{nm} X^{n} Y^{m}) = \sum c_{n0} X^{n}.$$
The units in both $A$ and $B$ are any
series with a nonzero constant
term, so
$A$ is SI in $B$.  Both $A$ and $B$ are local rings (in the
algebraic sense that they have a
unique maximal ideal), having only
one simple module given by evaluating at $X=Y=0$.
If we place the topology of
pointwise convergence of coefficients on
$A$ and $B$, then both algebras
are $Q$-algebras, since the
set of series with zero constant
term is clearly closed.
\vskip\baselineskip
{\vbox{\centerline {\S 2 Spectral
Invariance of $M_{n}({\Ainf})$ }
\vskip\baselineskip
\par
If $\Cal A$ is an algebra, let $M_{n}(\Cal A)$ denote the
$n$ by $n$ matrices over $\Cal A$.
{\theorem{Theorem 2.1.\ \ }}{\it  Let
$\varphi \colon \Ainf \rightarrow B$ be
an algebra homomorphism with dense image,
and assume that $B$ is a $Q$-algebra.  Then
$\Ainf$ is  SI in
$B$ iff $M_n(\Ainf)$ is  SI in $M_n(B)$.
}}
{\demo}
We prove the forward direction.   Assume
that $A$ is SI in $B$.
We will show that condition
$(\romannumeral3)$ of Theorem 1.4
holds for the pair $M_n(\Ainf)$,
$M_n(B)$.
We first need a lemma.
{\theorem{Lemma 2.2.\ \ }}{\it Let
$\Cal A$ be any algebra, $V$ a simple
$M_n(\Cal A)$-module.  Then $V$ is a direct sum
$$V = \bigoplus_1^n V_{0}$$
where $V_{0}$ is a simple $\Cal A$-module.  }
{\demo } Simple exercise using matrix units.\qed
\par
Let $V$ be a simple
$M_n(\Ainf)$-module.  By the lemma, $V=\bigoplus_1^n V_{0}$
where $V_{0}$ is a simple $\Ainf$-module.  Since $A$ is
SI in $B$ and $B$ is a $Q$-algebra,
we can extend $V_{0}$ to a $B$-module $W_{0}$.
Then $W=\bigoplus
_1^n W_{0}$ is an $M_n(B)$-module
in a natural way,  containing
$V$.
Hence $M_n(\Ainf)$
is SI in $M_n(B)$.
\qed
\par
Since $(\romannumeral3)\Rightarrow (\romannumeral1)$
in Theorem 1.4 did not require $B$ to be a $Q$-algebra,
we did not need to check that $M_{n}(B)$ is a $Q$-algebra
in Theorem 2.1.  It is true
by {\cite{[12, Lemma 2.1]}} that if $B$ is good, then
$M_{n}(B)$ is good.  However, I do not know if this is
true with ``good'' replaced by ``$Q$-algebra''.
{\theorem {Corollary 2.3.\ \ }}
{\it If $A$ is a dense Fr\'echet
subalgebra of a C*-algebra $B$,
then Theorem 2.1 is true
with ``spectral invariant'' replaced by
``closed under the holomorphic functional
calculus''.}
\par
{\bf Remark 2.4.\ \ } B. Gramsch has an unpublished
proof of Theorem 2.1 which
works if $B$ is a good topological
algebra with jointly continuous multiplication.
This proof was also noticed by
Bost in {\cite{[3, Proposition A.2.2]}},
and is very similar to the proof of {\cite{[7, Proposition 1.8.4]}}.
The proof fits roughly into
the following two steps.  Step 1 was first
proven by Swan in {\cite{[12, Lemma 1.1]}}.
\par
{\bf Step 1.\ \  } Assume both
$A$ and $B$ are unital (the nonunital
case follows easily from the
unital case).  For simplicity, we assume
$n=2$.  Since $B$ is a good topological algebra with jointly
continuous multiplication, we can find a
neighborhood $V$ of the identity in $M_{n}(B)$ such that
$$\left(\matrix {b_{11} & b_{12} \cr b_{21} & b_{22}}\right)\in V$$
implies $b_{11}$, $b_{22}$, and $b_{22}-b_{21}b^{-1}_{11}b_{12}$
are all invertible in $B$ (note we have used the fact that
inversion is continuous here).  If $A$ is
SI in $B$ and
$$[a_{ij}] = \left(\matrix {a_{11} & a_{12} \cr
a_{21} & a_{22} } \right)
\in V\cap M_{n}(A),$$
then $a_{11}$, $a_{22}$, $a_{22}- a_{21} a_{11}^{-1} a_{12}$ are
all invertible in $A$.  So
$$
\left(
\matrix {a_{11}^{-1} & 0 \cr 0 &
(a_{22}- a_{21}a_{11}^{-1}a_{12})^{-1}}
\right)
\left(
\matrix {1 & -a_{12}(a_{22}-a_{21} a_{11}^{-1} a_{12})^{-1}
\cr
0 & 1 } \right)\left(
\matrix {1 &  0 \cr -a_{21}a_{11}^{-1} & 1 }\right)
$$
is easily seen to be an inverse for $[a_{ij}]$.
(To construct this inverse, we have simply done elementary
row and column operations to put the matrix into diagonal
form, then cleared the diagonal.)
\par
{\bf Step 2.\ \  } Next we note
that if we can find such a neighborhood
$V$ as in step 1, then $M_{n}(A)$ is SI in $M_{n}(B)$.
For let $a\in M_{n}(A)\cap {\rm inv}M_{n}(B)$.  By density, there
are $a_{1}, a_{2} \in M_{n}(A)$ such that $a_{1}a\in V$,
$aa_{2}\in V$.  It follows that $a$ is invertible in $M_{n}(A)$.
This is the argument for the
last statement of {\cite{[8, Lemma 5.3]}},
which is also given in {\cite{[3, Proposition A.2.2]}}.
%
\par
{\bf Question 2.5.\ \ } Do there exist $Q$-algebras
 without continuous inversion ?
The field without continuous inversion
described in {\cite{[13, IX.2]}} may be such an example,
however I have
not been able to fill in the details of this
example, and have heard that others have had the same problem.
\vskip\baselineskip
\centerline{ \S 3 Examples }
\vskip\baselineskip
\par
{\bf Example 3.1.\ \ } We give an example of a non SI dense
Banach subalgebra of a C*-algebra
with an open group of invertible
elements.
\par
Let $B=C[-1,1]$.  Let $A$ be the set of elements of $B$ which
extend to analytic functions
on the open disc $\{\, z\in \C\, \,|\,\, |z|<1\, \}$,
continuous on the boundary.  The sup norm over the whole disc
makes $A$ a Banach algebra.  The subalgebra
$A$ is not SI since $i+z$ is invertible
in $B$.  For a simple $A$-module not extending to any $B$-module,
let $V=\C$, $fc = f(-i)c$.
\vskip\baselineskip
\par
We use Theorem 1.4 to show that some dense subalgebras
are {\it not} spectral invariant.
\par
{\bf Example 3.2.\ \ } Let $\Z$ act on
  $\T$ by an irrational rotation $\alpha$.
Let $B$ be the C*-crossed product $C^{*}(\Z, C(\T), \alpha)$,
and let $A=C_{c}^{\infty}(\Z \times \T)$,
the convolution algebra of $C^{\infty}$ functions with
compact support on $\Z \times \T$.  The
 representation of $A$ on $C^{\infty}(\T)$
given by
$$(Ff)(z) = \sum_{n\in \Z}
F(n,z)e^{-n} f(e^{-\iota \alpha n}z), \eqno (3.3)$$
where $F\in A$, $f\in C^{\infty}(\T)$ and  $z\in \T$,
is not contained in any $B$-module, since the character
$n\mapsto e^{-n}$ does not give a character of $C^{*}(\Z)$.
We show that (3.3) defines an algebraically irreducible
representation of $A$.  If $f\in C^{\infty}(\T)$ and
$f$ is not identically zero, then $|f| \in C^{\infty}(\T)f
\subseteq Af$.  So $Af$ contains a nonnegative
and nonzero function,
and all of its translates by multiples of $\alpha$.
Since $\T$ is compact, the sum of finitely many of these translates
never vanishes on $\T$.
Call this sum $g$.  Then $C^{\infty}(\T)\subseteq C^{\infty}(\T)g
\subseteq Af$.
Thus the representation is algebraically irreducible.
By Theorem 1.4,
$A$ is not spectral invariant in $B$. We remark that if the action
of $\Z$ on a manifold $M$ is {\it proper}, then the convolution
algebra $C_{c}^{\infty}(\Z \times M)$ {\it is} SI in the C*-crossed
product {\cite{[1, Appendix]}}.
\par
{\bf Example 3.4.\ \ } Let $G$ be the real $ax+b$ group,
defined as ${\bf R}^{2}$ with multiplication
$(a,b)(c,d) = (a+c, e^{a}d + b)$.  We define a dense
subalgebra of $C^{*}(G)$ consisting of ``exponentially decreasing''
Schwartz
functions, which we show is not spectral
invariant.
Define $\omega(a, b) = 1+e^{|a|}+|b| +|e^{-a}b|$.
Then $\omega$ is a submultiplicative function (or weight)
on $G$, and we define the $\omega$-rapidly
vanishing differentiable functions on $G$ to be those differentiable
functions $\psi$ from $G$ to $\C$ satisfying
$$
\int_{G} \omega^{d}(g) \, |D\psi (g) |\, dg < \infty
\eqno (3.5)$$
for each differential operator $D$ on $G$ and each natural number
$d$.  Here $dg$
denotes the left Haar measure on $G$.
This space ${\Cal S}^{\omega}(G)$
is then a dense $m$-convex Fr\'echet *-subalgebra of
$C^{*}(G)$ {\cite{[11, \S 1]}}.
Left Haar measure on $G$ is given by
$d\mu(a, b) = e^{-a} da db$, where $(a, b) \in \R^{2}$ and
$dadb$ is Lebesgue measure on $\R^{2}$.
The character
$$ \psi \mapsto %
 \int_{\R^{2}} \psi(a, b) da db
=\int e^{a} \psi(a, b) d\mu(a, b) \eqno (3.6) $$
defines a simple ${\Cal S}^{\omega}(G)$-module
not contained in any $C^{*}(G)$-module,
since the character $(a, b) \mapsto e^{a}$ of $G$ is not
unitary.
So ${\Cal S}^{\omega}(G)$ is not SI in $C^{*}(G)$.
\vskip\baselineskip
\par
We conclude with two examples of what happens in Theorem 1.4
$(i)\Rightarrow
(iii)$ when the condition that $B$ is a $Q$-algebra is dropped.
\par
{\bf Example 3.7.\ \  } We give an example of a pair of topological
algebras such that $A$ is SI and dense
in $B$, $\invB$ is not open, and there is
a simple $A$-module not contained in any $B$-module.
Let $B=\C[X, Y]$, $A= \C [Z, Y]$, the sets of complex polynomials
in $X, Y$ and $Z, Y$ respectively.
Identify
$A$ with a subalgebra of $B$ by mapping
$Z$ to $X(1-Y)$.
\par
We give $B$ the  following locally convex Hausdorff
topology.  Let $\C[[X, Y]]$
be the algebra of power series in two variables.
If $n, m \in \N$, we define ideals $I_{n, m}=
<X^{n}, Y^{m}>$
in $\C[[X, Y]]$.  Then we say that $f_{k}\rightarrow f_{0}$
in $\C[[X, Y]]$
if and only if $f_{k}-f_{0}$ is eventually in
every $I_{n, m}$.  Let
$B\subseteq \C[[X, Y]]$ have the relative topology.
With this topology on $B$,
we see that $X(1-Y)(1+Y +\dots Y^{n})$ converges
to $X$ as $n\rightarrow \infty$, so $X\in
{\overline A}^{B}$ and $A$ is dense in $B$.
The units in both
algebras are the constants, so $A$ is SI in $B$,
and $\invB$ is not open.
Evaluation at $Z$, $Y=1$ gives a simple $A$-module
contained in no $B$-module.
\par
{\bf Example 3.8.\ \  } We investigate further
what happens in Theorem 1.4 if we drop
the hypothesis that $\invB$ is open.
It is still true (in the implication $(\romannumeral1) \Rightarrow
(\romannumeral3)$) that if an extension exists it can be taken
simple.  However,
we show by example that extensions may always exist, but
we may not be able
to take our extension to be
continuous, or have dense inclusion.  In fact, in
our example $V$ is one-dimensional, but any extension
$W$ is discontinuous and infinite dimensional.
\par {\bf 1.\ \  }
Let $B=s(\bf N)$, the sequences of
(possibly unbounded) complex numbers
with pointwise multiplication and the topology of
pointwise convergence.  This space
is a complete, unital,
locally convex, $m$-convex Fr\'echet algebra over $\C$
(see  {\cite{[9, Example 3.7]}}).
An element of $B$ is invertible iff it is never zero.  One
easily finds a sequence of non-invertible elements converging
to $1$, so ${\rm inv}  B$
is not open.  Hence we know by  {\cite{[14]}} that $B$ will
have discontinuous infinite dimensional simple modules.
\par {\bf 2.\ \  }
We characterize all simple $B$-modules, or
equivalently all maximal ideals.
Let $\Cal U$ denote the
set of ultrafilters on $\bf N$ (see {\cite{[4]}}).
For any $ F \in \Cal U$ we define an ideal in $B$ by
$$I_{ F} = \{\, b \in B\,\, |\quad \exists S\in { F}
(b_n = 0 \Longleftrightarrow n\in S) \quad \}
\eqno (3.9) $$
One easily checks that this is a maximal ideal, and that every
maximal ideal arises in this way.  For each $p\in \bf N$, we have
a maximal ideal $I_p
= \{\, b\in B \,\, |\,\, b_p=0 \, \}$.  Such a maximal ideal
corresponds to the principal ultrafilter ${ F}_{p}
= \{\,S\subseteq \bf N\,\, |
\,\, p\in S\, \}$.  There are many
non-principal ultrafilters on $\bf N$
(see {\cite{[4, Theorem 7.1]}}).
\par
Any non-principal ultrafilter $ F$ contains all cofinite sets
so that $I_{ F}$ contains all finitely supported sequences
in $B$.  Thus $I_{ F}$ is a dense maximal ideal in $B$ and
$B/I_{ F}$ is a discontinuous simple $B$-module.
One can find elements of $B/I_{ F}$ which are not scalar
multiples of the identity, so that $B/I_{ F}$ is a proper
field extension of the algebraically closed field $\bf C$.
Hence $B/I_{ F}$ is infinite dimensional over $\bf C$.
\par {\bf 3.\ \  }
Now we define our dense subalgebra $A$ to be the algebra of
all finitely supported sequences $s_{\scriptscriptstyle f}(\bf N)$
with unit adjoined.  Let $V$ be the trivial $A$-module,
defined by $V=\bf C$ with action $(s + \lambda )z = \lambda z$,
for
$s\in
s_{\scriptscriptstyle f}
(\bf N)$,
$\lambda \in \bf C$, and
$z\in V$.  Clearly $s_{\scriptscriptstyle f}(\bf N)$
is the kernel of this representation, so the kernel of any
irreducible representation of $B$ containing $V$ is a maximal
ideal containing
$s_{\scriptscriptstyle f}(\bf N)$.  Hence any extension of $V$
to a simple $B$-module $W$ must be of the form $W=B/I_{ F}$
for some non-principal ultrafilter $ F$.  Hence there is
no extension of $V$ to a continuous simple $B$-module.
%
\par {\bf 4.\ \  }
We note, however, that every simple $A$-module does
have an extension to a $B$-module (so that $(\romannumeral 1)
\Rightarrow (\romannumeral 3)$ does not fail).
To show this, we classify
all simple $A$-modules by
noting that $s_{\scriptscriptstyle f}(\bf N)$ is
dense and
SI{} in the C{*}-algebra $c_0(\bf N)$ of sequences tending to
zero at infinity with the sup norm topology.  By Theorem 1.4, all
simple $A$-modules extend to simple $\tilde c_0(\bf N)$-
modules.  Thus the only simple $A$-modules are the
point evaluations and the trivial
module $V$ defined in {\bf 3}, all of which
extend to $B$-modules.
\vskip\baselineskip
\centerline{References}
\vskip\baselineskip
\frenchspacing
\item{[1]} P. Baum, A. Connes,
{\it Chern Character for Discrete Groups},
Collection: A fete of topology, Academic press,
Orlando, Florida,   1988, pp. 163-232.
\frenchspacing
\item{[2]}  B. Blackadar, $K$-theory for operator
algebras, Springer-Verlag, New York, 1986.
\frenchspacing
\item{[3]} J.B. Bost, {\it  Principe D'Oka, K-Theorie
et Systems Dynamiques Non-commutative},
Invent. Math, {\bf 101} (1990), 261-333.
\frenchspacing
\item{[4]} W.W. Comfort, S. Negrepontis,
The Thoery of Ultrafilters, Springer-Verlag,
Berlin/Heidleberg, 1974, pp.142-144.
\frenchspacing
\item{[5]} A. Connes,
{\it An Analogue of the Thom Isomophism for
Crossed Products of a C*-algebra by an Action of $\R$},
Adv. in Math.  {\bf 39} (1981), 31--55.
\frenchspacing
\item{[6]}   J. Dixmier,  C* algebras,
North-Holland Publishing Co., Amsterdam New York
Oxford,1982.
%
%
%
\frenchspacing
\item{[7]}   T. Palmer,
Banach Algebras and the
General Theory of *-Algebras,
Cambridge University Press Encyclopedia
of Mathematics Series, in preparation.
\frenchspacing
\item{[8]} B. Gramsch, {\it
Relative Inversion in der St\"orungstheorie
von Operatoren und $\Psi$-Algebren},
Math Ann  {\bf 209} (1984), 27--71.
\frenchspacing
\item{[9]} E. Micheal,  {\it Locally multiplicatively
convex topological algebras},   Mem. Amer. Math. Soc.,
{\bf 11} ( 1952 ).
\frenchspacing
\item{[10]}  N.C. Phillips,  {\it
$K$-theory for Fr\'echet algebras}, Intl. Jour.
Math. {\bf 2 (1)} (1991), 77-129.
\frenchspacing
\item{[11]}  L. Schweitzer, {\it
Dense $m$-convex Fr\'echet Subalgebras of Operator
Crossed Products by Lie Groups},  preliminary version
(1990).
\frenchspacing
\item{[12]} R. G. Swan, {\it Topological Examples
of Projective Modules},  Trans. Amer. Math. Soc.{\bf 230} (1977),
201-234.
\frenchspacing
\item{[13]}    L.   Waelbroeck,
Topological Vector Spaces
and Algebras,  Springer-Verlag,  Berlin Heidelberg
New York, 1971.
\frenchspacing
\item{[14]} W. Zelazko, {\it  On maximal ideals in commutative
$m$-convex algebras},   Studia Math {\bf T.LVIII} (1976),
291-298.
\end